\documentclass{article}
\usepackage{amsrefs}
\usepackage{amssymb}
\usepackage{amsmath}
\usepackage{amsthm}
\usepackage{graphicx}

\theoremstyle{definition}
\newtheorem{definition}{Definition}[section]

\newtheorem{remark}[definition]{Remark}
\theoremstyle{plain}

\newtheorem{proposition}[definition]{Proposition}
\newtheorem{theorem}[definition]{Theorem}

\DeclareMathOperator{\Div}{Div}
\DeclareMathOperator{\ord}{ord}

\author{Marc Coppens\footnote{Katholieke Hogeschool Kempen, Departement IBW,
Kleinhoefstraat 4, B-2440 Geel, Belgium; K.U.Leuven, Department of
Mathematics, Celestijnenlaan 200B, B-3001 Leuven, Belgium; email:
marc.coppens@khk.be} and Filip Cools\footnote{K.U.Leuven, Department of
Mathematics, Celestijnenlaan 200B, B-3001 Leuven, Belgium; email:
filip.cools@wis.kuleuven.be}}
\title{Linear pencils on graphs and on real curves}
\date{}

\begin{document}
\maketitle \noindent

\section{Introduction}

Let $X$ be a general smooth irreducible complex curve of genus $g$
and assume $r$ and $d$ are nonnegative integers such that $\rho
=g-(r+1)(g-d+r)=0$. From Brill-Noether Theory, it is known that $X$
has finitely many linear systems $g^r_d$ and that the natural
associated $0$-dimensional scheme $W^r_d$ parameterizing those
linear systems is reduced and has $\lambda
=g!\prod^r_{i=0}\frac{i!}{(g-d+r+i)!}$ points.

In \cite{ref1}, a new method for considering Brill-Noether problems
is introduced using linear systems on graphs. In particular, one
obtains a metric graph $\Gamma$ of genus $g$ having exactly
$\lambda$ linear systems $g^r_d$. In a degeneration from smooth
curves of genus $g$ to a singular curve of genus $g$ with dual graph
equal to $\Gamma$, the results from \cite{ref2} imply that the
generic curve $X$ has finitely many linear systems $g^r_d$. However,
it is not clear whether all linear systems $g^r_d$ on $X$ have
different specializations on $\Gamma$ (see \cite[Conjecture
1.5]{ref1}). In this note we solve this problem for the case $r=1$.

In \cite{ref4}, the specialization from curves to graphs is studied
for the situation of real curves. In particular, in that situation,
the graph $\Gamma$ has a real structure compatible to the real
structure on the curves and real linear systems on the curves
specialize to real linear systems on the graph. Using suited lengths
for the edges of the graph $\Gamma$ mentioned before, it is possible
to introduce two real structures on $\Gamma$, including the trivial
one. Using the previous results, we obtain the existence of real
curves of genus $g$ having exactly $\lambda$ linear systems $g_d^1$
all of them being real (using the trivial structure) and having less
than $\lambda$ real linear systems $g_d^1$ (using the non-trivial
structure) in case $g\geq 6$. In this second case, the number
$\lambda '$ of real linear systems $g_d^1$ can be computed and it turns out
that $\lambda '/\lambda$ becomes $0$ for $g\rightarrow \infty$.

The existence of real curves having $\lambda$ real linear systems
$g^1_d$ also follows from degenerations to real rational cuspidal
curves using \cite{ref5} and it is proved using the theory of limit
linear systems in \cite{ref6}. In general, one has the following
problem: determine all values $\lambda ' \equiv \lambda \pmod{2}$
such that there exists a real curve $X$ having exactly $\lambda$
linear systems $g^1_d$ with exactly $\lambda '$ of them being real.
The non-trivial real structure shows that $\lambda '$ can become
much smaller than $\lambda$. It should be noted that from the result
in \cite{ref7}, it follows that one cannot obtain the existence of
real curves with $\lambda '\neq \lambda$ using degenerations to real
cuspidal rational curves. In case $g$ is not of the type $2^n-2$,
then $\lambda$ is even and the smallest possible value of $\lambda
'$ could be $0$. Using the intersection of a general real cubic
surface and a general real quadric containing no real lines, it is
easy to find real curves of genus $4$ having no real $g^1_3$. In
\cite{ref8}, it is proved that there exist real curves of genus $8$
having no real linear system $g^1_5$. At the moment, those seem to
be the only known cases with $\lambda '=0$.

\section{Linear pencils on metric graphs}

Let $G$ be a finite connected graph, possibly with multiple edges,
but without loops. Denote by $V(G)$ the set of vertices of $G$ and
by $E(G)$ the set of edges of $G$. Assume that a real positive
number $\ell(e)$ is assigned to each edge $e$. If we identify each
edge $e\in E(G)$ with a line segment $[0,\ell(e)]$, we get a compact
connected metric space $\Gamma$, which is the metric graph
associated to $G$. The genus $g$ of $\Gamma$ is defined as
$|E(G)|-|V(G)|+1$, i.e. the first Betti number of $G$, where $G$ is
any underlying graph.

A divisor $D=n_1v_1+\ldots+n_tv_t$ on $\Gamma$ is an element of the
free abelian group $\Div(\Gamma)$ on the points of the metric graph
$\Gamma$. The degree of $D$ is the sum $n_1+\ldots+n_t$ of the
coefficients and $D$ is effective if each coefficient $n_i$ is
nonnegative. Often, the terminology of chips is used when one deals
with divisors on graphs, e.g. we say that the divisor $D$ has $k$
chips at $v$ if the coefficient of $v$ in $D$ is equal to $k$.

Let $\psi:\Gamma\to\mathbb{R}$ be a continuous map such that $\psi$
is piecewise linear with finitely many pieces and with integer
slopes on each edge of $\Gamma$. The divisor $\text{div}(\psi)$ of
$\psi$ is equal to $\sum_{v\in\Gamma}\,\ord_v(\psi)$, where
$\ord_v(\psi)$ is the sum of the incoming slopes of $\psi$ at
$v\in\Gamma$. We say that two divisors $D$ and $D'$ are equivalent
if and only if $D'-D=\text{div}(\psi)$ for some continuous map
$\psi$ on $\Gamma$, and in this case, we write $D\sim D'$. The set
$|D|$ of all effective divisors $D'$ equivalent with $D$ is the
linear system corresponding to $D$. The rank $r(D)$ of $|D|$ or $D$
is defined as follows. If $|D|=\emptyset$, then $r(D)=-1$, otherwise
$r(D)\geq r$ if and only if $|D-E|\neq\emptyset$ for each effective
divisor $E$ on $\Gamma$ of degree $r$. For a more elaborate
introduction to linear systems on metric graphs (e.g. the notion of
the $v$-reduced divisor of $|D|$), we refer to \cite{ref1}. \\

Linear systems on metric graphs seem to obey analogous results as
linear systems on algebraic curves, e.g. Riemann-Roch Theorem. In
\cite{ref2}, Baker conjectures that also the following counterpart
of the Brill-Noether Theorem holds and he proves the first part.

\begin{theorem}[Brill-Noether Theorem for metric graphs]
Fix $g,r,d\geq 0$ and set $\rho=g-(r+1)(g-d+r)$.
\begin{enumerate}
\item[(1)] If $\rho\geq 0$, then every metric graph $\Gamma$ of genus $g$ has
a divisor $D$ of rank $r(D)=r$ and degree $\deg(D)\leq d$.
\item[(2)] If $\rho<0$, then there exists a metric graph $\Gamma$
of genus $g$ for which there does not exist a divisor $D$ of rank
$r(D)=r$ and degree $\deg(D)\leq d$.
\end{enumerate}
\end{theorem}

The second part was proved in \cite{ref1} using the following metric
graph. Let $G_g$ be the graph consisting of $g$ loops that are
connected via $4$-valent vertices $v_1,\ldots,v_{g-1}$ and such that
in addition, the first and the last loop also contain one $2$-valent
node each, respectively $v_0$ and $v_g$. Assume that $\Gamma$ is the
corresponding metric graph where the lengths $\ell_i$ and $m_i$ of
the two edges from $v_{i-1}$ to $v_i$ are generic for each
$i\in\{1,\ldots,g\}$ (it suffices to take $\ell_i=\ell\geq 2g-2$ and
$m_i=1$ for all $i$). In \cite{ref1} is proved that such metric
graphs are Brill-Noether general.

\begin{figure}[h]
\centering
\includegraphics[width=0.7\textwidth]{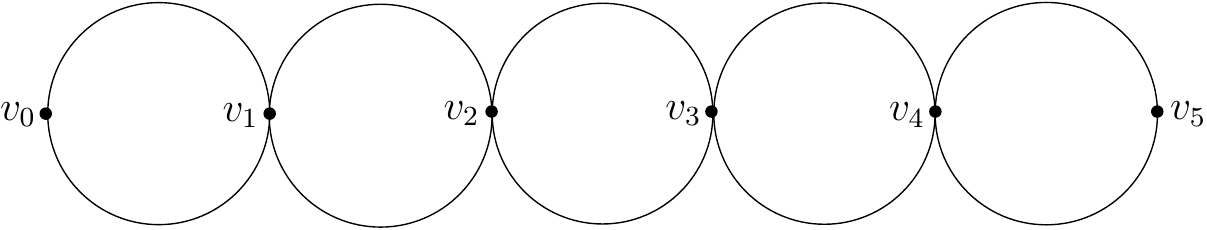}
\caption{the graph $G_5$}
\end{figure}

We will now focuss on the case of linear pencils on $\Gamma$,
thus $r=1$. If $g$ is even and $d=\frac{g}{2}+1$ (so $\rho=0$),
it is proven in \cite{ref1} that there is a bijection between
lattice paths $p=(p_0,\ldots,p_g)$ in $\mathbb{Z}$ satisfying
$p_0=p_g=1$, $p_i\geq 1$ and $p_i-p_{i-1}=\pm 1$ for all
$i\in \{1,\ldots,g\}$ and linear pencils of degree $d$ on
$\Gamma$ as follows. If $p=(p_0,\ldots,p_g)$ is a path satisfying
the conditions, let $D_p$ be the divisor on $\Gamma$ with one chip
in $v_0$, one (extra) chip on the unique point $w_i$ of the $i$th
loop satisfying $p_{i-1}v_{i-1}+w_i\sim p_iv_i$ if $p_i-p_{i-1}=1$
and no (extra) chips on the $i$th loop if $p_i-p_{i-1}=-1$.
Note that $D_p$ is $v_0$-reduced. The bijection maps $p$ to the
linear system $|D_p|$. So the number of linear pencils on $\Gamma$
is equal to the Catalan number $\lambda=\frac{1}{d} {2d-2 \choose d-1}$.

\begin{proposition} \label{propdouble}
Let $D_p$ be a divisor on $\Gamma$ corresponding to a lattice path
$p=(p_0,\ldots,p_g)$. Then the linear system $|2D_p|$ has rank equal to two.
\end{proposition}
\begin{proof}
Since it is clear that $|2D_p|$ has rank at least two, it suffices
to prove that $|2D_p-2v_0-v_g|=\emptyset$. Let $Q_i$ be the divisor
attained from $2D_p-2v_0$ by moving as many chips as possible from
the first $i$ loops of $\Gamma$ to $v_i$, for each
$i\in\{0,\ldots,g\}$. Denote by $q_i$ the number of chips of $Q_i$
in $v_i$. We claim that $q_i=p_i-1$ for each $i$ and we will prove
this by induction on $i$.

For $i=0$, the divisor $Q_0=2D_p-2v_0$ and thus $q_0=0=p_0-1$ since
$D_p$ has exactly one chip in $v_0$. Now assume the claim holds for
$i$, hence $q_i=p_i-1$. If $p_{i+1}-p_i=-1$, the divisor $Q_{i+1}$
is attained from $Q_i$ by moving the $q_i$ chips in $v_i$ as much as
possible to $v_{i+1}$. Since the edge lengths are general, only
$q_i-1=p_i-2=p_{i+1}-1$ chips end up in $v_{i+1}$. If
$p_{i+1}-p_i=1$, the divisor $Q_i$ has $q_i$ chips in $v_i$ and $2$
chips at the point $w_{i+1}$ of the $(i+1)$th loop satisfying
$p_iv_i+w_{i+1}\sim p_{i+1}v_{i+1}$, so only
$(q_i+2)-1=p_i=p_{i+1}-1$ chips end up in $v_{i+1}$. In both cases,
we have $q_{i+1}=p_{i+1}-1$.

The statement of the proposition now follows from the claim for $i=g$,
since the $v_g$-reduced divisor $Q_g$ in $|2D-2v_0|$ has $q_g=p_g-1=0$
chips in $v_g$.
\end{proof}

Assume that the metric graph $\Gamma$ has integer edge lengths (thus
$\ell_i,m_i\in\mathbb{Z}$). Let $G$ be the graph corresponding to
$\Gamma$ with $\ell(e)=1$ for each $e\in E(G)$, so $G$ is a
refinement of $G_g$. Let $R$ be a complete discrete valuation ring
with field of fractions $Q$ and algebraically closed residue field
$k$. Then \cite[Appendix B]{ref2} implies that there exists a
regular arithmetic surface $\mathfrak{X}\to R$ whose generic fiber
is a smooth curve $X\to Q$ and whose special fiber $\mathfrak{X}_k$
has dual graph equal to $G$. Following \cite[Section 2]{ref2}, the
specialization map $$\alpha:\Div(X)\to\Div(G)$$ sends effective (resp.
principal) divisors to effective (resp. principal) divisors. Hereby,
$\Div(G)$ is the subgroup of $\Div(\Gamma)$ consisting of divisors
supported at $V(G)$. Let $\overline{Q}$ be the algebraic closure of
$Q$. The map $\alpha$ can be extended to a map
$$\tau_*:\Div(X(\overline{Q}))\cong \Div(X_{\overline{Q}})\to
\Div(\Gamma)$$ satisfying $$r(\tau_*(E))\geq \dim(|E|)$$ for each
$E\in \Div(X_{\overline{Q}})$. This inequality implies that also
$X_{\overline{Q}}$ is Brill-Noether general, but it is not clear
wether each divisor on $\Gamma$ of rank $r$ and degree $d$ lifts to
a divisor of rank $r$ and degree $d$ on $X_{\overline{Q}}$ (see
\cite[Conjecture 1.5]{ref1}). Here we solve this problem for the
case where $r=1$ and $\rho=0$.

\begin{theorem} \label{thmbijection}
If $g$ is even and $d=\frac{g}{2}+1$, each linear pencil $|D_p|$ on
$\Gamma$ lifts to a unique linear pencil $g_d^1$ on
$X_{\overline{Q}}$.
\end{theorem}
\begin{proof}
If $|E|$ would be a multiple $g_d^1$ on $X_{\overline{Q}}$, then it
follows from the description of the tangent map to $W_d^1$ and the
base point free pencil trick (see \cite[Chapter IV, Section 4 and
Chapter III, Section 3]{ref3}) that $\dim(|2E|)\geq 3$. This would
imply that $|E|$ specializes to a linear pencil $|D_p|$ on $\Gamma$
satisfying $$r(2D_p)=r(\tau_*(2E))\geq \dim(|2E|)\geq 3,$$ which is
in contradiction with Proposition \ref{propdouble}. So the scheme
$W_d^1$ of $X_{\overline{Q}}$ consists of exactly
$\lambda=\frac{1}{d}{2d-2 \choose d-1}$ points.

Assume two different linear pencils $|E_1|$ and $|E_2|$ on
$X_{\overline{Q}}$ specialize to the same linear pencil $|D_p$ on
$\Gamma$. Then $\dim(|E_1+E_2|)\geq 3$ and $r(2D_p)=2$, but on the
other hand
$$r(2D_p)=r(\tau_*(E_1+E_2))\geq \dim(|E_1+E_2|),$$ so we have a
contradiction.

Since the graph $\Gamma$ has exactly $\lambda$ pencils $g_d^1$, it
follows that $\tau_*$ induces a bijection between linear pencils
$g_d^1$ on $X_{\overline{Q}}$ and on $\Gamma$.
\end{proof}

Assume $\ell_i=\ell$ and $m_i=1$ for all $i\in\{1,\ldots,g\}$ and
identify the two edges of $\Gamma$ between $v_{i-1}$ and $v_i$ with
intervals $I_i=[0,1]$ and $J_i=[0,\ell]$. In this case, there is a
non-trivial involution $\sigma$ on $\Gamma$ sending $x\in I_i$ to
$1-x\in I_{g-i}$ and $x\in J_i$ to $\ell-x\in J_{g-i}$. Note that the
only fixed point of this involution is $v_{g/2}$.

\begin{proposition} \label{propsigma}
Let $p=(p_0,\ldots,p_g)$ be a lattice path corresponding to a
divisor $D_p$ of rank one on $\Gamma$. Denote by $\sigma(p)$ the
lattice path $(p_g,\ldots,p_0)$.
\begin{enumerate}
\item[(i)] $D_p\sim \sigma(D_{\sigma(p)})$.
\item[(ii)] The linear system $|D_p|$ is invariant under $\sigma$
if and only if $p=\sigma(p)$.
\end{enumerate}
\end{proposition}
\begin{proof}
Define $f_i$ for all $i\in\{0,\ldots,g\}$ inductively by $f_0=0$ and
$$f_i=\begin{cases} f_{i-1}+p_{i-1} & \text{if $p_i-p_{i-1}=1$} \\ f_{i-1}+p_i & \text{if $p_i-p_{i-1}=-1$} \end{cases}$$
Let $f:\Gamma\to\mathbb{R}$ be the function defined by
$$f(x)=\begin{cases} f_{i-1}+p_{i-1}.x & \text{if $x\in I_i$}  \\
f_{i-1} & \text{if $x\in [0,\ell-p_{i-1}]\subset J_i$}  \\
f_{i-1}+ x - (\ell-p_{i-1}) & \text{if $x\in [\ell-p_{i-1},\ell]\subset J_i$}  \end{cases}$$
if $p_i-p_{i-1}=1$ and
$$f(x)=\begin{cases} f_{i-1}+p_i.x & \text{if $x\in I_i$}
\\ f_{i-1}+x & \text{if $x\in [0,p_i]\subset J_i$}
\\ f_{i-1}+ p_i & \text{if $x\in [p_i,\ell]\subset J_i$}  \end{cases}$$
if $p_i-p_{i-1}=-1$. It is clear that $f$ is a piecewise linear with $f(v_i)=f_i$
and one can see that $\text{div}(f)=\sigma(D_{\sigma(p)})-D$, hence $D_p\sim \sigma(D_{\sigma(p)})$.

If the linear system $|D_p|$ is invariant under $\sigma$,
then $\sigma(D_{\sigma(p)})\in |D_p|$ implies that $D_{\sigma(p)}\in |D_p|$.
Since $D_p$ and $D_{\sigma(p)}$ are both $v_0$-reduced divisors in $|D_p|$, they must be equal, hence $p=\sigma(p)$.

Now assume $p=\sigma(p)$, hence $\sigma(D_p)\in|D_p|$. Let $E$ be an
arbitrary divisor in $|D_p|$, so there exists a piecewise linear
function $\psi:\Gamma\to\mathbb{R}$ such that
$\text{div}(\psi)=E-D_p$. If we apply the involution $\sigma$ to
this equation, we get that
$\text{div}(\psi\circ\sigma)=\sigma(\text{div}(\psi))=\sigma(E-D_p)=\sigma(E)-\sigma(D_p)$,
hence $\sigma(E)\sim\sigma(D_p)\sim D_p$, so $\sigma(E)\in |D_p|$ and $|D_p|$ is invariant under the
involution.
\end{proof}

This implies that the number $\lambda'$ of linear pencils on
$\Gamma$ that are invariant under $\sigma$ is equal to the number
of symmetric lattice paths. This number can be computed as follows.
Using Andr\'e's reflection method (see \cite{ref9}), one can see that
the number of symmetric lattice paths $p=(p_0,\ldots,p_g)$ with $p_{g/2}=m$
(where $m\in \{1,\ldots,d\}$ with $m\equiv d \mod 2$) is equal to
$$\frac{m}{d}{d \choose \frac{d-m}{2}}={ d-1 \choose \frac{d-m}{2}}-{d-1 \choose \frac{d-m}{2}-1},$$
hence $\lambda'$ is the sum of these numbers and equals the central binomial
coefficient ${d-1 \choose \lceil \frac{d-1}{2}\rceil}$. In the following table,
the values of $\lambda$ and $\lambda'$ are mentioned for small values of $d$.
Note that the limit of the quotient $\lambda'/\lambda$ goes to zero if $d\to\infty$.

\begin{center}
\begin{tabular}{|c||c|c|c|c|c|c|c|c|c|}
\hline $d$&2&3&4&5&6&7&8&9&10 \\ \hline $g$&2&4&6&8&10&12&14&16&18 \\ \hline $\lambda$&1&2&5&14&42&132&429&1430&4862 \\ \hline $\lambda'$&1&2&3&6&10&20&35&70&126 \\ \hline
\end{tabular}
\end{center}

\section{Linear pencils on real curves}

In \cite{ref4}, it is proved that real linear systems on curves
degenerate to real linear systems on graphs. We are going to explain
that those arguments also imply that the degeneration respects
complex conjugation of linear systems.

Let $R$ be a complete discrete valuation ring defined over
$\mathbb{R}$ having residue field $\mathbb{R}$ and let $\mathfrak{X}
\rightarrow R$ be a regular arithmetic surface such that the generic
fiber $X$ is smooth and geometrically irreducible. Let $Q$ be the
quotient field of $R$, then $Q_{\mathbb{C}}=Q\otimes
_{\mathbb{R}}\mathbb{C}$ is a field (see \cite[Lemma 5.1]{ref4}).
Let $R_{\mathbb{C}}$ be the complete valuation ring extending $R$
in $Q_{\mathbb{C}}$ and let $\mathfrak{X}_{\mathbb{C}} \rightarrow
R_{\mathbb{C}}$ be obtained by base change. We assume the special
fiber $X_{0,\mathbb{C}}$ (which is defined over $\mathbb{R}$) is
strongly semistable. This means $X_{0,\mathbb{C}}$ is reduced, all
singular points of $X_{0,\mathbb{C}}$ are nodes and its irreducible
components are smooth. Associated to $X_{0,\mathbb{C}}$, there is a
dual graph $G$ and this dual graph has a real structure defined by
the real structure on $X_{0,\mathbb{C}}$.

Let $D$ be a divisor on the generic fiber $X_{\mathbb{C}}$ and let
$\overline{D}$ be the conjugated divisor. This means
$\overline{D}=\sigma (D)$, where $\sigma$ is the involution on
$X_{\mathbb{C}}$ induced by the field extension $Q\subset
Q_{\mathbb{C}}$. Considering $X_{\mathbb{C}}\subset
\mathfrak{X}_{\mathbb{C}}$, one has closures $cl(D)$ and
$cl(\overline{D})$ and one has $cl(\overline{D})=\overline{cl(D)}$
(in the right hand side of this equality, the conjugation is induced
by the field extension $\mathbb{R}\subset \mathbb{C}$ because
$\mathfrak{X}$ is defined over $\mathbb{R}$).

Let $\overline{Q}$ be the algebraic closure of $Q$ and let
$\overline{Q}^r$ be the real closure. We obtain a field extension
$\overline{Q}^r\subset \overline{Q}$ of degree $2$. On
$\overline{X}=X\times _Q \overline{Q}$, it defines an involution
$\sigma$ and if $D$ is a divisor on $\overline{X}$, then we write
$\overline{D}$ to denote $\sigma (D)$. The divisor $D+\overline{D}$
is defined over $\overline{Q}^r$. In particular, there exists a
finite extension $Q\subset K$ with $K\subset \overline{Q}^r$ such
that $D+\overline{D}$ is defined over $K$ and $D$ is defined over
$K\otimes _{\mathbb{R}} \mathbb{C}=K_{\mathbb{C}}$. Let $R_K$ (resp.
$R_{K,\mathbb{C}}$) be the extension of $R$ in $K$ (resp.
$K_{\mathbb{C}}$) and consider the base extensions
$\mathfrak{X}\times _R R_K$ and $\mathfrak{X}\times _R
R_{K,\mathbb{C}}$. It is proved in \cite{ref4} that there exists a
family $\mathfrak{X}_K \rightarrow R_K$ defined over $\mathbb{R}$
such that by base change we obtain a desingularization
$\mathfrak{X}_{K,\mathbb{C}}\rightarrow R_{K,\mathbb{C}}$ of
$\mathfrak{X}\times _R R_{K,\mathbb{C}}$, the special fiber
$X_{0,K,\mathbb{C}}$ is strongly semistable and its graph $G_K$ with
the associated real stucture has a natural weighted structure such
that the associated metric graph is equal to $\Gamma$ with its real
structure. The associated divisor on $G_K$ defines a divisor
$\tau_*(D)$ on $\Gamma$ and one has
$\tau_*(\overline{D})=\overline{\tau_*(D)}$ on $\Gamma$. In case
$D=\overline{D}$, these arguments are explained in full detail in
\cite{ref4} showing that $\tau_*(D)$ is a real divisor on $\Gamma$
in that case.

We consider the metric graph $\Gamma$ used in the proof of
Proposition \ref{propsigma}. It has two real structures: the trivial
one and the one defined by $\sigma$. From \cite[Proposition
5.9]{ref4}, it follows that both real structures can be obtained
from a degeneration with $X_{0,\mathbb{C}}$ a totally degenerated
curve. This means that all components of $X_{0,\mathbb{C}}$ have
genus $0$, in particular
$g(X_{0,\mathbb{C}})=g(\overline{X})=g(\Gamma)=g$. Now assume such a
degeneration is chosen.

\begin{proposition}
The specialization induces a bijection between the set of real
linear systems $g_d^1$ on $\overline{X}$ and the set of real linear
systems $g_d^1$ on $\Gamma$.
\end{proposition}
\begin{proof}
Theorem \ref{thmbijection} implies that the specialization map
induces a bijection between the sets of all linear systems $g_d^1$.
From \cite{ref4}, it follows that under the specialization a real
$g_d^1$ on $\overline{X}$ corresponds to a real $g_d^1$ on $\Gamma$.
In case $g$ is a non-real $g_d^1$ on $\overline{X}$, then $g$ and
$\overline{g}$ are two different linear systems on $\overline{X}$,
hence $\tau_*(g)\neq \tau_*(\overline{g})$. Using the equality
$\tau_*(\overline{g})=\overline{\tau_*(g)}$ mentioned above, we get
$\overline{\tau_*(g)}\neq \tau_*(g)$, hence $\tau_*(g)$ is non-real
on $\Gamma$.
\end{proof}

Using Proposition \ref{propsigma} in the case the real structure on
$\Gamma$ is induced by $\sigma$ and using a general specialization
of the generic fiber of $\mathfrak{X}\to R$ to a complex curve $X$
of genus $g$ defined over $\mathbb{R}$, we obtain the following.

\begin{theorem}
There exist smooth complex curves $X$ defined over $\mathbb{R}$ of
genus $g=2d-2$ having exactly $\lambda$ linear systems $g_d^1$ such
that all of them (resp. exactly $\lambda'$ of them) are real.
\end{theorem}

\begin{remark}
In case the real structure on $\Gamma$ is induced by $\sigma$, we
can say a little bit more about the structure of the real curve. We
can choose $X_{0,\mathbb{C}}$ such that the component $C_{g/2}$
corresponding to $v_{g/2}$ has a non-empty real locus or an empty
real locus. The other components are obtained from pairs of disjoint
complex conjugated copies of $\mathbb{P}^1$. In the first case, the
deformation is a smooth real curve of genus $g$ having exactly one
connected component in the real locus. In the second case, we obtain
a real curve with an empty real locus. In the case of real curves
with empty real locus, there exist two types of linear systems
invariant under complex conjugation distinguished by the fact
whether or not they contain a real divisor. Unfortunately, using our
methods, we are not able to distinguish between those possibilities.
\end{remark}

\end{document}